\RequirePackage{ifpdf}
\ifpdf 
\documentclass[pdftex]{sigma}
\else
\documentclass{sigma}
\fi

\def\non{\nonumber}
\def\ep{\epsilon}
\def\lar{\longrightarrow}

\def\tila{\tilde{a}}
\def\tilb{\tilde{b}}
\def\tilc{\tilde{c}}
\def\tR{\tilde{R}}
\def\hR{\widehat{R}}
\def\tF{\tilde{F}}
\def\hF{\widehat{F}}
\def\tG{\tilde{G}}
\def\hG{\widehat{G}}
\def\tW{\tilde{W}}
\def\tPsi{\tilde{\Psi}}
\def\ttT{\tilde{{\mathbb T}}}

\def\pl{l^{+}}
\def\pr{r^{+}}

\begin{document}


\renewcommand{\PaperNumber}{049}

\FirstPageHeading

\ShortArticleName{Free Field Approach to Solutions of the Quantum
Knizhnik--Zamolodchikov Equations}

\ArticleName{Free Field Approach to Solutions\\ of the Quantum
Knizhnik--Zamolodchikov Equations}

\Author{Kazunori KUROKI~$^\dag$ and Atsushi NAKAYASHIKI~$^\ddag$}

\AuthorNameForHeading{K. Kuroki and A. Nakayashiki}

\Address{$^\dag$~Department of Mathematics, Kyushu University,
Hakozaki 6-10-1, Fukuoka 812-8581, Japan}
\EmailD{\href{mailto:ma306012@math.kyushu-u.ac.jp}{ma306012@math.kyushu-u.ac.jp}}

\Address{$^\ddag$~Department of Mathematics, Kyushu University,
Ropponmatsu 4-2-1, Fukuoka 810-8560, Japan$\!$}
\EmailD{\href{mailto:6vertex@math.kyushu-u.ac.jp}{6vertex@math.kyushu-u.ac.jp}}

\ArticleDates{Received February 18, 2008, in f\/inal form May 27,
2008; Published online June 03, 2008}

\Abstract{Solutions of the qKZ equation associated with the
quantum af\/f\/ine algebra $U_q(\widehat{sl}_2)$ and its two
dimensional evaluation
 representation are studied.
The integral formulae derived from the free f\/ield realization of
intertwining operators of $q$-Wakimoto modules are shown to
coincide with those of Tarasov and Varchenko.}


\Classification{81R50; 20G42; 17B69}

\section{Introduction}
In 1992 I.~Frenkel and N.~Reshetikhin \cite{FR} had developed the
theory of intertwining operators for quantum af\/f\/ine algebras
and had shown that the matrix elements of intertwiners satisfy the
quantized Knizhnik--Zamolodchikov (qKZ) equations.

The theory of intertwiners and qKZ equations was successfully
applied to the study of solvable lattice models
 \cite{JM} (and references therein).
As to the study of solutions of the qKZ equations, bases are
constructed by Tarasov and Varchenko \cite{TV} in the form of
multi-dimensional hypergeometric integrals in the case of
$U_q(\widehat{sl}_2)$. However solutions of the qKZ equations for
other quantum af\/f\/ine algebras are not well studied \cite{TV0}.

The method of free f\/ields is ef\/fective to compute correlation
functions in conformal f\/ield theo\-ry (CFT) \cite{ATY}, in
particular, solutions to the Knizhnik--Zamolodchikov (KZ)
equa\-tions~\mbox{\cite{SV1,SV2}}. A similar role is expected for
those of quantum af\/f\/ine algebras. Unfortunately it is
dif\/f\/icult to say that this expectation is well realized, as we
shall explain below.

Free f\/ield realizations of quantum af\/f\/ine algebras are
constructed by Frenkel and Jing \cite{FJ} for level one integrable
representations of ADE type algebras and by Matsuo~\cite{M},
Shiraishi~\cite{S} and Abada et al.~\cite{ABG} for representations
with arbitrary level of $U_q(\widehat{sl}_2)$. The latter results
are extended to $U_q(\widehat{sl}_N)$ in \cite{AOS}. Free f\/ield
realizations of intertwiners are constructed based on these
representations in the case of $U_q(\widehat{sl}_2)$
\cite{JM,KQS,K1,M,BW}.

The simplicity of the Frenkel--Jing realizations makes it possible
not only to compute matrix elements but also traces of
intertwining operators \cite{JM}, which are special solutions to
the qKZ equations. The case of $q$-Wakimoto modules with an
arbitrary level becomes
 more complex and the detailed study of the solutions of the qKZ equations
making use of it is not well developed. In \cite{M} Matsuo derived
his integral formulae \cite{M0} from the formulae obtained by the
free f\/ield calculation in the simplest case of one integration
variable. However it is not known in general whether the integral
formulae derived from the free f\/ield realizations recover those
of \cite{M0,V,TV}\footnote{See Note~1 in the end of the paper.}.

The aim of this paper is to study this problem in the case of the
qKZ equation with the value in the tensor product of two
dimensional irreducible representations of $U_q(\widehat{sl}_2)$.
More general cases will be studied in a subsequent paper.

There are mainly two reasons why the comparison of two formulae is
dif\/f\/icult. One is that the formulae derived from the free
f\/ield calculations contain more integration variables than in
Tarasov--Varchenko's (TV) formulae. This means that one has to
carry out some integrals explicitly to compare two formulae. The
second reason is that the formulae from free f\/ields  contains a
certain sum. This stems from the fact that the current and
screening operators are written as a sum which is absent in the
non-quantum case. Since TV formulae have a similar structure to
those for the solutions of the KZ equation \cite{SV1,SV2}, one
needs to sum up certain terms explicitly for the comparison of two
formulae. We carry out such calculations in the case we mentioned.

The plan of this paper is as follows. In Section~2 the
construction of the hypergeometric solutions of the qKZ equation
due to Tarasov and Varchenko is reviewed. The free f\/ield
construction of intertwining operators is reviewed in Section~3.
In Section~4 the formulae for the highest to highest matrix
elements of some operators are calculated. The main theorem is
also stated in this section. The transformation of the formulae
from free f\/ields to Tarasov--Varchenko's formulae is described
in Section~5. In Section~6 the proof of the main theorem is given.
Remaining problems are discussed in Section~7. The appendix
contains the list of the operator product expansions which is
necessary to derive the integral formula.

\section[Tarasov-Varchenko's formula]{Tarasov--Varchenko's formula}

Let $V^{(1)}={\mathbf C}v_0\oplus{\mathbf C}v_1$ be a
two-dimensional irreducible representation of the algebra
$U_q(sl_2)$, and $R(z)\in {\rm End}(V^{(1)\otimes 2})$ be a
trigonometric quantum $R$-matrix given by
\begin{gather*}
R(z)\left(v_\ep\otimes v_\ep\right) = v_\ep\otimes v_\ep, \non
\\
R(z)\left(v_0\otimes v_1\right)
 =
\frac{1-z}{1-q^2z}q v_0\otimes v_1 + \frac{1-q^2}{1-q^2z}
v_1\otimes v_0, \non
\\
R(z)\left(v_1\otimes v_0\right)
 =
\frac{1-q^2}{1-q^2z}z v_0\otimes v_1 + \frac{1-z}{1-q^2z}q
v_1\otimes v_0, \non
\end{gather*}
Let $p$ be a complex number such that $|p|<1$ and $T_j$ denote the
multiplicative $p$-shift operator of $z_j$,
\begin{gather*}
T_jf(z_1,\dots ,z_n)=f(z_1,\dots ,pz_j,\dots ,z_n).
\end{gather*}
The qKZ equation for the $V^{(1)\otimes n}$-valued function
$\Psi(z_1,\dots ,z_n)$ is
\begin{gather}
T_j\Psi=R_{j,j-1}(pz_j/z_{j-1})\cdots R_{j,1}(pz_j/z_{1})
\kappa^{\frac{1-h_j}{2}} R_{j,n}(z_j/z_n)\cdots
R_{j,j+1}(z_j/z_{j+1})\Psi, \label{qKZ1}
\end{gather}
where $R_{ij}(z)$ signif\/ies that $R(z)$ acts on the $i$-th and
$j$-th components, $\kappa$ is a complex parameter,
$\kappa^{\frac{1-h_j}{2}}$ acts on the $j$-th component as
\begin{gather*}
\kappa^{\frac{1-h_j}{2}}v_\ep=\kappa^\ep v_\ep.
\end{gather*}
Let us brief\/ly recall the construction of the hypergeometric
solutions \cite{TV,T} of the equation~\eqref{qKZ1}. In the
remaining part of the paper we assume $|q|<1$. We set
\begin{gather*}
(z)_\infty=(z;p)_\infty, \qquad
(z;p)_\infty=\prod_{j=0}^\infty(1-p^jz), \qquad
\theta(z)=(z)_\infty(pz^{-1})_\infty(p)_\infty.
\end{gather*}

Let $n$ and $l$ be non-negative integers satisfying $l\leq n$. For
a sequence $(\ep)=(\ep_1,\dots ,\ep_n)\in\{0,1\}^n$ satisfying
$\sharp\{i|\ep_i=1\}=l$ let
\begin{gather*}
w_{(\ep)}(t,z)= \prod_{a<b}   \frac{t_a-t_b}{q^{-2}t_a-t_b}
\sum_{1\leq a_1,\dots ,a_l\leq l \atop a_i\neq a_j (i\neq j)}
\prod_{i=1}^l \left( \frac{t_{a_i}}{t_{a_i}-q^{-1}z_{k_i}}
\prod_{j<k_i}\frac{q^{-1}t_{a_i}-z_j}{t_{a_i}-q^{-1}z_j}
\prod_{i<j}\frac{q^{-2}t_{a_i}-t_{a_j}}{t_{a_i}-t_{a_j}}\! \right)
,
\end{gather*}
where $\{i|\ep_i=1\}=\{k_1<\cdots<k_l\}$.

The elliptic hypergeometric space ${\cal F}_{ell}$ is the space of
functions $W(t,z)=W(t_1,\dots ,t_l,$ $z_1, \dots ,z_n)$ of the
form
\begin{gather*}
W=Y(z)\Theta(t,z)
\frac{1}{\prod_{j=1}^n\prod_{a=1}^l\theta(qt_a/z_j)} \prod_{1\leq
a<b\leq l}\frac{\theta(t_a/t_b)}{\theta(q^{-2}t_a/t_b)}
\end{gather*}
satisfying the following conditions
\begin{enumerate}\itemsep=0pt
\item[(i)] $Y(z)$ is meromorphic on $({\mathbb C}^\ast)^n$ in
$z_1$,\dots ,$z_n$, where ${\mathbb C}^\ast={\mathbb
C}\backslash\{0\}$; \item[(ii)] $\Theta(t,z)$ is holomorphic on
$({\mathbb C}^\ast)^{n+l}$ in $t_1$,\dots ,$z_n$ and symmetric in
$t_1$,\dots ,$t_l$; \item[(iii)] $T^t_aW/W=\kappa q^{n-2l+4a-2}$,
$ T^z_jW/W=q^{-l}$, where $T^t_a W=W(t_1,\dots ,pt_a,\dots t_l,z)$
and $T^z_j W=W(t,z_1,\dots ,pz_j,\dots z_n)$.
\end{enumerate}

Def\/ine the phase function $\Phi(t,z)$ by
\begin{gather*}
\Phi(t,z)=\prod_{i=1}^n\prod_{a=1}^l
\frac{(qt_a/z_i)_\infty}{(q^{-1}t_a/z_i)_\infty}
\prod_{a<b}\frac{(q^{-2}t_a/t_b)_\infty}{(q^{2}t_a/t_b)_\infty}.
\end{gather*}

For $W\in {\cal F}_{ell}$ let
\begin{gather*}
I(w_{(\ep)},W)=
\int_{\ttT^l}\prod_{a=1}^l\frac{dt_a}{t_a}\Phi(t,z)
w_{(\ep)}(t,z)W(t,z),
\end{gather*}
where $\ttT^l$ is a suitable deformation of the torus
\begin{gather*}
{\mathbb T}^l=\{(t_1,\dots ,t_l)|\, |t_i|=1, 1\leq i\leq l\},
\end{gather*}
specif\/ied as follows~\cite{TV}.

Notice that the integrand has simple poles at
\begin{gather*}
t_a/z_j = (p^sq^{-1})^{\pm1}, \qquad s\geq 0, \quad 1\leq a\leq l,
\quad 1\leq j\leq n, \non
\\
t_a/t_b = (p^sq^2)^{\pm1}, \qquad s\geq 0, \quad 1\leq a<b\leq l.
\end{gather*}
The contour for the integration variable $t_a$ is a simple closed
curve which rounds the origin in the counterclockwise direction
and separates the following two sets,
\begin{gather*}
\{p^sq^{-1}z_j, p^sq^2t_b\,|\, s\geq 0,\, 1\leq j\leq n,\,
a<b\,\}, \non
\\
\{p^{-s}qz_j,p^{-s}q^{-2}t_b\,|\, s\geq 0,\, 1\leq j\leq n,\,
a<b\,\}.
\end{gather*}
Then
\begin{gather*}
\Psi_W=\sum_{(\ep)}I(w_{(\ep)},W)v_{\ep_1}\otimes \cdots \otimes
v_{\ep_n},
\end{gather*}
is a solution of the qKZ equation \eqref{qKZ1} for any $W\in {\cal
F}_{ell}$.

\section[Free field realizations]{Free f\/ield realizations}

In this section we review the free f\/ield construction of the
representation of the quantum af\/f\/ine algebra
$U_q(\widehat{sl_2})$ of level $k$ and intertwining operators. We
mainly follow the notation in \cite{KQS}. We set
\begin{gather*}
[x]=\frac{q^x-q^{-x}}{q-q^{-1}}.
\end{gather*}
Let $k$ be a complex number and
$\{a_n,b_n,c_n,\tila_0,\tilb_0,\tilc_0,Q_a,Q_b,Q_c|n\in {\mathbb
Z}\backslash\{0\}\}$ satisfy
\begin{alignat*}{3}
& [a_n,a_m]=\delta_{m+n,0} \frac{[2n][(k+2)n]}{n}, \qquad&&
[\tila_0,Q_a]=2(k+2),& \non
\\
& [b_n,b_m]= -\delta_{m+n,0} \frac{[2n]^2}{n}, \qquad &&
[\tilb_0,Q_b]=-4,& \non
\\
& [c_n,c_m]= \delta_{m+n,0} \frac{[2n]^2}{n}, \qquad &&
[\tilc_0,Q_c]=4, & \non
\end{alignat*}

Other combinations of elements are supposed to commute. Set
\begin{gather*}
N_\pm={\mathbb C}[a_n,b_n,c_n|\pm n>0].
\end{gather*}
Then the Fock module $F_{r,s}$ is def\/ined to be the free
$N_{-}$-module of rank one generated by the vector which
satisf\/ies
\begin{gather*}
N_{+}|r,s\rangle  = 0, \qquad \tila_0|r,s \rangle = r|r,s\rangle ,
\qquad \tilb_0|r,s \rangle = -2s|r,s\rangle , \qquad \tilc_0|r,s
\rangle = 2s|r,s\rangle .
\end{gather*}
We set
\begin{gather*}
F_r=\oplus_{s\in {\mathbb Z}}F_{r,s}.
\end{gather*}
A representation of the quantum af\/f\/ine algebra
$U_q(\widehat{sl_2})$ is constructed on $F_r$ for any $r\in
{\mathbb C}$ in~\cite{S}.

The right Fock module $F_{r,s}^\dag$ and $F_r^\dag$ are similarly
def\/ined using the vector $\langle r,s|$ satisfying the
conditions
\begin{gather*}
\langle r,s|N_{-} = 0, \qquad \langle r,s|\tila_0 = r\langle r,s|,
\qquad \langle r,s|\tilb_0 = -2s\langle r,s|, \qquad \langle
r,s|\tilc_0 = 2s\langle r,s|.
\end{gather*}

\begin{remark} We change the def\/inition of $|r,s\rangle $ in \cite{KQS}.
Namely we use
\begin{gather*}
|r,s\rangle
=\exp\left(\frac{r}{2(k+2)}Q_a+s\frac{Q_b+Q_c}{2}\right)|0,0\rangle.
\end{gather*}
\end{remark}

Let us introduce f\/ield operators which are relevant to our
purpose. For $x=a,b,c$ let
\begin{gather*}
x(L;M,N|z:\alpha) = -\sum_{n\neq
0}\frac{[Ln]x_n}{[Mn][Nn]}z^{-n}q^{|n|\alpha}+
\frac{L\tilde{x}_0}{MN}\log\,z+\frac{L}{MN}Q_x, \non
\\
x(N|z:\alpha) = x(L;L,N|z:\alpha) = -\sum_{n\neq
0}\frac{x_n}{[Nn]}z^{-n}q^{|n|\alpha}+
\frac{\tilde{x}_0}{N}\log\,z+\frac{1}{N}Q_x.
\end{gather*}
The normal ordering is def\/ined by specifying $N_{+}$, $\tila_0$,
$\tilb_0$, $\tilc_0$ as annihilation operators and $N_{-}$, $Q_a$,
$Q_b$, $Q_c$ as creation operators. With this notation let us
def\/ine the operators
\begin{gather*}
J^{-}(z): F_{r,s} \lar  F_{r,s+1}, \qquad \phi^{(l)}_m(z): F_{r,s}
\lar  F_{r+l,s+l-m}, \qquad S(z): F_{r,s} \lar  F_{r-2,s-1},
\end{gather*}
by
\begin{gather*}
J^{-}(z)=\frac{1}{(q-q^{-1})z}\left(J^{-}_{+}(z)-J^{-}_{-}(z)\right),
\non
\\
J^{-}_{\mu}(z)= :\exp\left(
a^{(\mu)}\left(q^{-2}z;-\frac{k+2}{2}\right)
+b\big(2|q^{(\mu-1)(k+2)}z;-1\big)
+c\big(2|q^{(\mu-1)(k+1)-1}z;0\big) \right):, \non
\\
a^{(\mu)}\left(q^{-2}z;-\frac{k+2}{2}\right)= \mu\left\{
(q-q^{-1})\sum_{n=1}^\infty a_{\mu n}z^{-\mu
n}q^{(2\mu-\frac{k+2}{2})n}+\tila_0\log\, q\right\}, \non
\\
S(z)=\frac{-1}{(q-q^{-1})z}\left(S_{+}(z)-S_{-}(z)\right), \non
\\
S_{\ep}(z)=:\exp\left( -a\left(k+2|q^{-2}z;-\frac{k+2}{2}\right)
-b\big(2|q^{-k-2}z;-1\big) -c\big(2|q^{-k-2+\epsilon}z;0\big)
\right):,
\\
\phi^{(l)}_l(z)= :\exp\left(
a\left(l;2,k+2|q^kz;\frac{k+2}{2}\right) \right):, \non
\\
\phi^{(l)}_{l-r}(z)= \frac{1}{[r]!}
\oint\prod_{j=1}^r\frac{du_j}{2\pi i}
\left[\cdots\left[\left[\phi^{(l)}_l(z),J^{-}(u_1)\right]_{q^l},J^{-}(u_2)
\right]_{q^{l-2}}, \dots, J^{-}(u_r)\right]_{q^{l-2r+2}}, \non
\end{gather*}
where
\begin{gather*}
[r]!=[r][r-1] \cdots [1], \qquad [X,Y]_q=XY-qYX,
\end{gather*}
and the integral in $\phi^{(l)}_{l-r}(z)$ signif\/ies to take the
coef\/f\/icient of $(u_1\cdots u_r)^{-1}$.

The operator $J^{-}(z)$ is a generating function of a part of
generators of the Drinfeld realization $U_q(\widehat{sl_2})$ at
level $k$. While the operators $\phi^{(l)}_m(z)$ are conjectured
to determine the intertwining operator for $U_q(\widehat{sl_2})$
modules \cite{KQS,M}
\begin{gather*}
\phi^{(l)}(z):W_r\lar W_{r+l}\otimes V^{(l)}_z, \qquad
\phi^{(l)}(z)=\sum_{m=0}^l\phi^{(l)}_m(z)\otimes v^{(l)}_m,
\end{gather*}
where $W_r$ is a certain submodule of $F_r$ specif\/ied as a
kernel of a certain operator, called $q$-Wakimoto module
\cite{M,K1,K2,K3,ABG}, $V^{(l)}$ is the irreducible representation
of $U_q(sl_2)$ with spin~$l/2$ and $V^{(l)}_z$ is the evaluation
representation of $U_q(\widehat{sl_2})$ on $V^{(l)}$.

In this paper we exclusively consider the case $l=1$ and set
\begin{gather*}
\phi_{+}(z)=\phi^{(1)}_0, \qquad \phi_{-}(z)=\phi^{(1)}_1, \qquad
v_{0}=v^{(1)}_0, \qquad v_{1}=v^{(1)}_1.
\end{gather*}

The operator $S(z)$ commutes with $U_q(\widehat{sl_2})$ modulo
total dif\/fe\-ren\-ce. Here modulo total dif\/ference means
modulo functions of the form
\begin{gather*}
{}_{k+2} \partial _z
f(z):=\frac{f(q^{k+2}z)-f(q^{-(k+2)}z)}{(q-q^{-1})z}.
\end{gather*}

\begin{remark} The intertwining properties of $\phi^{(l)}(z)$
for $l\in {\mathbb Z}$ are not proved in \cite{KQS} as pointed out
in \cite{M}. However the fact that the matrix elements of
compositions of $\phi^{(l)}(z)$'s and $S(t)$'s satisfy the qKZ
equation modulo total dif\/ference can be proved in a similar way
to Proposition~6.1 in~\cite{M} using the result of Konno~\cite{K3}
(see~\eqref{qKZ2}).
\end{remark}

Let
\begin{gather*}
|m \rangle =|m,0\rangle \in F_{m,0}, \qquad \langle m|= \langle
m,0|\in F_{m,0}^\dag.
\end{gather*}
They become left and right highest weight vectors of
$U_q(\widehat{sl_2})$ with the weight $m\Lambda_1+(k-m)\Lambda_0$
respectively, where $\Lambda_0$, $\Lambda_1$ are fundamental
weights of $\widehat{sl}_2$. Consider
\begin{gather}
F(t,z)= \langle
m+n-2l|\phi^{(1)}(z_1)\cdots\phi^{(1)}(z_n)S(t_1)\cdots
S(t_l)|m\rangle \label{F}
\end{gather}
which is a function taking the value in $V^{(1)\otimes n}$. Let
\begin{gather*}
\Delta_j=\frac{j(j+2)}{4(k+2)}, \qquad s=\frac{1}{2(k+2)}.
\end{gather*}
Set
\begin{gather}
\hF=\left(\prod_{i=1}^nz_i^{\Delta_{m+n-2l+1-i}-\Delta_{m+n-2l-i}}\right)\,F
=\left(\prod_{i=1}^nz_i^{s(m+n-2l-i+\frac{3}{2})}\right)\,F,
\label{hF}
\end{gather}
Let the parameter $p$ be def\/ined from $k$ by
\begin{gather*}
p=q^{2(k+2)}.
\end{gather*}
We assume $|p|<1$ as before. Then the function $\hF$ satisf\/ies
the following qKZ equation modulo total dif\/ference of a function
\cite{M,IIJMNT,FR,KQS,K3}
\begin{gather}
T_j^z\hF= \hR_{j,j-1}(pz_j/z_{j-1}){\cdots} \hR_{j,1}(pz_j/z_{1})
q^{-(m+n/2-l+1)h_j} \hR_{j,n}(z_j/z_n){\cdots}
\hR_{j,j+1}(z_j/z_{j+1})\hF,\!\!\!\! \label{qKZ2}
\end{gather}
where
\begin{gather*}
\hR(z)=\rho(z)\tR(z), \qquad \tR(z)=C^{\otimes 2}R(z)C^{\otimes
2}, \non
\\
\rho(z)=q^{1/2}
\frac{(z^{-1};q^4)_\infty(q^4z^{-1};q^4)_\infty}{(q^2z^{-1};q^4)_\infty^2},
\qquad (z;x)_\infty=\prod_{i=0}^\infty(1-x^iz), \qquad
Cv_{\ep}=v_{1-\ep}.
\end{gather*}

\section{Integral formulae}

Def\/ine the components of $F(t,z)$ by
\begin{gather*}
F(t,z)=\sum_{(\nu)\in\{0,1\}^n}F^{(\nu)}(t,z)v_{(\nu)}, \qquad
v_{(\nu)}= v_{\nu_1}\otimes\cdots\otimes v_{\nu_n},
\end{gather*}
where $(\nu)=(\nu_1,\dots ,\nu_n)$. By the weight condition
$F^{(\nu)}(t,z)=0$ unless the condition
\begin{gather*}
\sharp\{i|\nu_i=0\}=l
\end{gather*}
is satisf\/ied. We assume this condition once for all. Notice that
\begin{gather*}
\phi_{+}(z)=\frac{1}{(q-q^{-1})}\oint\frac{du}{2\pi
iu}[\phi_{-}(z),J^{-}_{+}(u)- J^{-}_{-}(u)]_q, \non
\\
S(t)=\frac{-1}{(q-q^{-1})t}(S_{+}(t)-S_{-}(t)).
\end{gather*}
Let
\begin{gather*}
\{i|\nu_i=0\}=\{k_1<\cdots<k_l\},
\end{gather*}
and
\begin{gather*}
F^{(\nu)}_{(\ep)(\mu)}(t,z|u)= \langle
m+n-2l|\phi_{-}(z_1)\cdots[\phi_{-}(z_{k_1}),J^{-}_{\mu_1}(u_1)]_{q}\cdots
[\phi_{-}(z_{k_l}),J^{-}_{\mu_l}(u_l)]_{q} \cdots\phi_{-}(z_n)
\\
\phantom{F^{(\nu)}_{(\ep)(\mu)}(t,z|u)=}{}  \times
S_{\ep_1}(t_1)\cdots S_{\ep_l}(t_l)|m\rangle .
\end{gather*}
Then $F^{(\nu)}(t,z)$ can be written as
\begin{gather*}
F^{(\nu)}(t,z)=(-1)^l(q-q^{-1})^{-2l}\prod_{a=1}^lt_{a}^{-1}
\sum_{\ep_i,\mu_j}\prod_{i=1}^l(\ep_i\mu_i)
\int_{C^l}\prod_{j=1}^l\frac{du_j}{2\pi i u_j}
F^{(\nu)}_{(\ep)(\mu)}(t,z|u),
\end{gather*}
where $C^l$ is a suitable deformation of the torus ${\mathbb T}^l$
specif\/ied as follows.

The contour for the integration variable $u_i$ is a simple closed
curve rounding the origin in the counterclockwise direction such
that $q^{k+3}z_j$ $(1\leq j\leq n)$, $q^{-2}u_j$ $(i<j)$,
$q^{-\mu_i(k+2)}t_a\, (1\leq a\leq l)$ are inside and $q^{k+1}z_j$
$(1\leq j\leq n)$, $q^{2}u_j$ $(j<i)$ are outside.

By the operator product expansions (OPE) of the products of
$\phi_{-}(z)$, $J^{-}_{\mu}(u)$, $S_{\ep}(t)$ in the appendix, one
can compute the function $F^{(\nu)}_{(\ep)(\mu)}(t,z|u)$
explicitly. In order to write down the formula we need some
notation. Set
\begin{gather*}
\xi(z)=\frac{(pz^{-1};p,q^4)_\infty(pq^4z^{-1};p,q^4)_\infty}{(pq^2z^{-1};p,q^4)_\infty},
\qquad
(z;p,q)_\infty=\prod_{i=0}^\infty\prod_{j=0}^\infty(1-p^iq^jz).
\end{gather*}
Then
\begin{gather*}
F^{(\nu)}_{(\ep)(\mu)}(t,z|u)=f^{(\nu)}_{(\mu)}(t,z|u)\Phi(t,z)
G^{(\nu)}_{(\ep)(\mu)}(t,z|u),
\end{gather*}
where
\begin{gather*}
f^{(\nu)}_{(\mu)}(t,z|u) =(1-q^2)^l
q^{\sum_{i=1}^l(n+m-2l-k_i+i)\mu_i} \prod_{i=1}^n
\big(q^kz_i\big)^{s(m+n-l-i)} \prod_{i<j}\xi(z_i/z_j)
\\
\phantom{f^{(\nu)}_{(\mu)}(t,z|u)=}{} \times
\prod_{a=1}^l(q^{-2}t_a)^{4s(a-1)-2ms},
\\
G^{(\nu)}_{(\ep)(\mu)}(t,z|u) = \hG^{(\nu)}_{(\ep)(\mu)}(t,z|u)
\prod_{a<b}\frac{q^{\ep_b}t_b-q^{\ep_a}t_a}{t_b-q^{-2}t_a},
\\
\hG^{(\nu)}_{(\ep)(\mu)}(t,z|u) = \prod_{i=1}^l
\frac{u_i(z_{k_i}-q^{\mu_i-2-k}u_i)}
{(z_{k_i}-q^{-1-k}u_i)(u_i-q^{k+3}z_{k_i})}
\prod_{j=1}^l\prod_{i<k_j}\frac{z_i-q^{\mu_j-2-k}u_j}{z_i-q^{-1-k}u_j}
\\
\phantom{\hG^{(\nu)}_{(\ep)(\mu)}(t,z|u)=}{} \times
\prod_{j=1}^l\prod_{k_j<i}\frac{u_j-q^{k+2-\mu_j}z_i}{u_j-q^{k+3}z_i}
\prod_{i<j}\frac{u_i-q^{\mu_i-\mu_j}u_j}{u_i-q^{-2}u_j}
\prod_{i,a}\frac{u_i-q^{-\mu_i(k+1)-\ep_a}t_a}{u_i-q^{-\mu_i(k+2)}t_a}.
\end{gather*}
Let
\begin{gather*}
G^{(\nu)}_{(\mu)}(t,z)=\sum_{\ep_1,\dots ,\ep_n=\pm}
\prod_{j=1}^n\ep_j\int_{C^l}\prod_{j=1}^l\frac{du_j}{2\pi i u_j}
G^{(\nu)}_{(\ep)(\mu)}(t,z).
\end{gather*}
The main theorem in this paper is

\begin{theorem}
If $(\mu)\neq(-^l)=(-,\dots ,-)$, $G^{(\nu)}_{(\mu)}(t,z)=0$. For
$(\mu)=(-^l)$ we have
\begin{gather*}
G^{(\nu)}_{(-^l)}(t,z) =
q^{-2l+\frac{1}{2}l(l-1)-\sum_{i=1}^lk_i}(q-q^{-1})^l
w_{(-\nu)}(t,z),
\end{gather*}
where $(-\nu)=(1-\nu_1,\dots ,1-\nu_n)$.
\end{theorem}

It follows that $F^{(\nu)}(t,z)$ is given by
\begin{gather*}
F^{(\nu)}(t,z)=(-1)^l
q^{-(n+m+2-2l)l+ksn(m+n-l)-\frac{1}{2}ksn(n+1)+4sl(m-l+1)}
\\
\phantom{F^{(\nu)}(t,z)=}{} \times \prod_{i=1}^nz_i^{s(m+n-l-i)}
\prod_{i<j}\xi(z_i/z_j) \prod_{a=1}^lt_a^{2s(2a-2-m)-1} \Phi(t,z)
w_{(-\nu)}(t,z).
\end{gather*}

\section[Transformation to Tarasov-Varchenko's formulae]{Transformation to Tarasov--Varchenko's formulae}

We describe a transformation from $F$, which satisf\/ies
\eqref{qKZ2},
 to $\Psi$, which satisf\/ies \eqref{qKZ1}. The parame\-ter~$\kappa$ is also
determined as a function of $l$, $m$, $n$.

For a solution $\hG$ of \eqref{qKZ2} let
\begin{gather}
 \tG = \prod_{i=1}^nz_i^{-s(m+\frac{n}{2}-l+1)}
\left((\prod_{i=1}^nz_i)^{s/2}\prod_{i<j}\xi(z_i/z_j)\right)^{-1}
C^{\otimes n}\hG. \label{tG}
\end{gather}
One can easily verify that $\tG$ satisf\/ies \eqref{qKZ1} with
$\kappa=q^{2l-2-n-2m}$ using
\begin{gather*}
\frac{\xi(pz)}{\xi(z)}=\frac{(z^{-1};q^4)_\infty(q^4z^{-1};q^4)_\infty}{(q^2z^{-1};q^4)_\infty^2}.
\end{gather*}
Let $\hF$ be def\/ined by \eqref{hF} and $\tF$ by \eqref{tG}. Then
\begin{gather*}
\tF =
(-1)^lq^{-(n+m+2-2l)l+ksn(m+n-l)-\frac{1}{2}ksn(n+1)+4sl(m-l+1)}
\non
\\
\phantom{\tF =}{} \times \prod_{i=1}^nz_i^{s(m+n-2l-i)}
\prod_{a=1}^lt_a^{2s(2a-2-m)-1} \Phi(t,z) \sum
w_{(\nu)}(t,z)\otimes v_{(\nu)}.
\end{gather*}
For $W\in {\cal F}_{ell}$ let
\begin{gather}
\tW=W\left(\prod_{i=1}^nz_i^{s(m+n-2l-i)}
\prod_{a=1}^lt_a^{2s(2a-2-m)}\right)^{-1}. \label{tW}
\end{gather}
Then the condition (iii) for $W$ is equivalent to the following
conditions,
\begin{gather*}
T^t_a\tW/\tW=1, \qquad T^z_j\tW/\tW=q^{l-m-n+j}.
\end{gather*}

To sum up we have

\begin{proposition}
For any $W\in {\cal F}_{ell}$
\begin{gather*}
\tPsi_W=\int_{\ttT^l}\prod_{a=1}^l\frac{dt_a}{2\pi
i}\tF(t,z)\tW(t,z),
\end{gather*}
is a solution of the qKZ equation \eqref{qKZ1}, where $\tF$ is
defined by~\eqref{tG} with $\hF$ and $F$ being given in~\eqref{hF}
and~\eqref{F} and $\tW$ is defined by~\eqref{tW}.
\end{proposition}

\section{Proof of Theorem}

Let $A^\pm=\{j|\mu_j=\pm\}$. Suppose that the number of elements
in $A^\pm$ is $r^\pm$ and
 write $A^\pm=\{l^\pm_1<\cdots<l^\pm_{r^\pm}\}$.
Set $A=A^{-}$, $r=r^{-}$ and $l_i=l_i^{-}$for simplicity. Let
\begin{gather*}
I^{(\nu)}_{(\ep)(\mu)}(t,z)=
\int_{C^l}\prod_{j=1}^l\frac{du_j}{2\pi i
u_j}\hG^{(\nu)}_{(\ep)(\mu)}(t,z|u).
\end{gather*}

\begin{lemma}
We have
\begin{gather*}
I^{(\nu)}_{(\ep)(\mu)}(t,z) = (-q^{-2})^r(q-q^{-1})^r
\sum_{a_1,\dots ,a_r=1\atop a_i\neq a_j (i\neq j)}^l
\prod_{i=1}^r\delta_{\ep_{a_i},+}
\prod_{i=1}^r\frac{t_{a_i}}{z_{k_{l_i}}-qt_{a_i}}
\prod_{j=1}^r\prod_{i<k_{l_j}}\frac{z_i-q^{-1}t_{a_j}}{z_i-qt_{a_j}}
\non
\\
\phantom{I^{(\nu)}_{(\ep)(\mu)}(t,z)=}{}
\times\prod_{i=1}^r\prod_{a\neq a_i,\dots ,a_r}
\frac{t_{a_i}-q^{-1-\ep_a}t_a}{t_{a_i}-t_a}.
\end{gather*}
\end{lemma}

\begin{proof}
We f\/irst integrate in the variables $u_j$, $j\in A^{+}$ in the
order $u_{\pl_1},\dots ,u_{\pl_{\pr}}$. Let us consider the
integration in $u_{\pl_1}$. We denote the integration contour in
$u_i$ by $C_i$. The only singularity of the integrand outside
$C_{\pl_1}$ is $\infty$. Thus the integral is calculated by taking
residue at $\infty$. Since the integrand is of the form
\begin{gather*}
\frac{du_{\pl_1}}{u_{\pl_1}}H(u_{\pl_1}),
\end{gather*}
where $H(u)$ is holomorphic at $\infty$. Then
\begin{gather*}
\int_{C_{\pl_1}}\frac{du}{2\pi i u}H(u)=-\mathop{\rm
Res}_{u=\infty}\frac{du}{u}H(u)= \lim_{u\lar \infty}H(u).
\end{gather*}
In this way the integral in $u_{\pl_1}$ is calculated. After this
integration the integrand as a function of $u_{\pl_2}$ has a
similar structure. Therefore the integration with respect to
$u_{\pl_2}$ is carried out in a~similar way and so on. Finally we
get
\begin{gather*}
I^{(\nu)}_{(\ep)(\mu)}(t,z) = (-1)^{\pr} \mathop{\rm
Res}_{u_{\pl_{\pr}}=\infty} \cdots \mathop{\rm
Res}_{u_{\pl_{1}}=\infty} \hG^{(\nu)}_{(\ep)(\mu)}(t,z|u) \non
\\
\phantom{I^{(\nu)}_{(\ep)(\mu)}(t,z)}{} =
\int_{C^{n-\pr}}\prod_{j\in A}\frac{du_j}{2\pi i u_j} \prod_{i\in
A}\frac{-q^{-3-k}u_i}{z_{k_i}-q^{-1-k}u_i} \prod_{i<k_j \atop j\in
A}\frac{z_i-q^{-3-k}u_j}{z_i-q^{-1-k}u_j} \prod_{i<j\atop i,j\in
A}\frac{u_i-u_j}{u_i-q^{-2}u_j} \non
\\
\phantom{I^{(\nu)}_{(\ep)(\mu)}(t,z)=}{} \times \prod_{i\in A\atop
a}\frac{u_i-q^{k+1-\ep_a}t_a}{u_i-q^{k+2}t_a}. \non
\end{gather*}
Here $C^{n-\pr}$ is specif\/ied by similar conditions to $C^l$,
where $u_{\pl_i}$ $1\leq i\leq \pr$ are omitted. We denote the
right hand side of this equation other than
$\int_{C^{n-\pr}}\prod_{j}\frac{du_j}{2\pi i u_j}$ by
$I^{(\nu)+}_{(\ep)(\mu)}(t,z)$.

Next we integrate with respect to the remaining variables $u_j$,
$j\in A$ in the order $u_{l_r}$,\dots ,$u_{l_1}$. Let us consider
the integration with respect to $u_{l_r}$. The poles of the
integrand inside $C_{l_r}$ is $q^{k+2}t_a$, $a=1,\dots ,l$. Thus
we have
\begin{gather*}
\int_{C_{l_r}}\frac{du_{l_r}}{2\pi i
u_{l_r}}I^{(\nu)+}_{(\ep)(\mu)}(t,z) = \sum_{a_r=1}^l \mathop{\rm
Res}_{u_{l_r}=q^{k+2}t_{a_r}}I^{(\nu)+}_{(\ep)(\mu)}(t,z)
\\
\quad{}= \sum_{a_r=1}^l \prod_{i\in A \atop i\neq
l_r}\frac{-q^{-3-k}u_i}{z_{k_i}-q^{-1-k}u_i} \prod_{i<k_j \atop
j\in A\backslash\{l_r\}} \frac{z_i-q^{-3-k}u_j}{z_i-q^{-1-k}u_j}
\prod_{i<j<l_r\atop i,j\in A}\frac{u_i-u_j}{u_i-q^{-2}u_j}
\prod_{i\in A\backslash\{l_r\}\atop a}
\frac{u_i-q^{k+1-\ep_a}t_a}{u_i-q^{k+2}t_a}
\\
\qquad{}\times \frac{-q^{-1}t_{a_r}}{z_{k_{l_r}}-qt_{a_r}}
\prod_{i<k_{l_r}}\frac{z_i-q^{-1}t_{a_r}}{z_i-qt_{a_r}}
\prod_{i<k_{l_r}\atop i\in
A}\frac{u_i-q^{k+2}t_{a_r}}{u_i-q^{k}t_{a_r}} \prod_{a\neq
a_r}\frac{t_{a_r}-q^{-1-\ep_a}t_a}{t_{a_r}-t_a} (1-q^{-1-\ep_a})
\\
\quad{}= (1-q^{-2}) \sum_{a_r=1}^l \delta_{\ep_{a_r},+}
\prod_{i\in A \atop i\neq
l_r}\frac{-q^{-3-k}u_i}{z_{k_i}-q^{-1-k}u_i} \prod_{i<k_j \atop
j\in A\backslash\{l_r\}}\frac{z_i-q^{-3-k}u_j}{z_i-q^{-1-k}u_j}
\prod_{i<j<l_r\atop i,j\in A}\frac{u_i-u_j}{u_i-q^{-2}u_j}
\\
\qquad{}\times \left(\prod_{i\in A\backslash\{l_r\}\atop a\neq
a_r} \frac{u_i-q^{k+1-\ep_a}t_a}{u_i-q^{k+2}t_a}\right)
\frac{-q^{-1}t_{a_r}}{z_{k_{l_r}}-qt_{a_r}}
\prod_{i<k_{l_r}}\frac{z_i-q^{-1}t_{a_r}}{z_i-qt_{a_r}}
\prod_{a\neq a_r}\frac{t_{a_r}-q^{-1-\ep_a}t_a}{t_{a_r}-t_a}.
\end{gather*}
The last expression as a function of $u_{l_{r-1}}$ has a similar
form to $I^{(\nu)+}_{(\ep)(\mu)}(t,z)$ as a function of $u_{l_r}$.
Namely the poles inside $C_{l_{r-1}}$ are $q^{k+2}t_a$, $a\neq
a_r$. Thus the integral in $u_{l_{r-1}}$ is the sum of residues at
$q^{k+2}t_a$, $a\neq a_r$ and so on. Finally we get
\begin{gather*}
I^{(\nu)}_{(\ep)(\mu)}(t,z) = \sum_{a_r=1}^l\sum_{a_{r-1}=1\atop
a_{r-1}\neq a_r}^l\cdots \sum_{a_1\atop a_1\neq a_2,\dots ,a_r}^l
\mathop{\rm Res}_{u_{l_1}=q^{k+2}t_{a_1}} \cdots \mathop{\rm
Res}_{u_{l_r}=q^{k+2}t_{a_r}}I^{(\nu)+}_{(\ep)(\mu)}(t,z)
\\
\phantom{I^{(\nu)}_{(\ep)(\mu)}(t,z)}{} = (-q^{-2})^r(q-q^{-1})^r
\sum_{a_1,\dots ,a_r=1\atop a_i\neq a_j (i\neq j)}^l
\prod_{i=1}^r\delta_{\ep_{a_i},+}
\prod_{i=1}^r\frac{t_{a_i}}{z_{k_{l_i}}-qt_{a_i}}
\prod_{j=1}^r\prod_{i<k_{l_j}}\frac{z_i-q^{-1}t_{a_j}}{z_i-qt_{a_j}}
\\
\phantom{I^{(\nu)}_{(\ep)(\mu)}(t,z)=}{} \times
\prod_{i=1}^r\prod_{a\neq a_i,\dots ,a_r}
\frac{t_{a_i}-q^{-1-\ep_a}t_a}{t_{a_i}-t_a}.
\end{gather*}
Thus the lemma is proved.
\end{proof}

Recall that
\begin{gather*}
G^{(\nu)}_{(\ep)(\mu)}(t,z) = I^{(\nu)}_{(\ep)(\mu)}(t,z)
\prod_{a<b}\frac{q^{\ep_b}t_b-q^{\ep_a}t_a}{t_b-q^{-2}t_a}.
\end{gather*}
For a set $\{a_1,\dots ,a_r\}$ let $\{b_1,\dots ,b_{\pr}\}$ be
def\/ined by
\begin{gather*}
\{1,\dots ,l\}=\{a_i\}\cup\{b_i\}.
\end{gather*}
Then
\begin{gather*}
(-q^2)^r(q-q^{-1})^{-r}G_{(\mu)}^{(\nu)}(t,z)
\\
= \sum_{a_1,\dots ,a_r=1 \atop a_i\neq a_j (i\neq j)}^l\!
\prod_{i=1}^r\! \frac{t_{a_i}}{z_{k_{l_i}}-qt_{a_i}} \prod_{i=1}^r
\prod_{j<k_{l_i}}\frac{z_j-q^{-1}t_{a_i}}{z_j-qt_{a_i}}
\prod_{i>j}\frac{t_{a_i}-q^{-2}t_{a_j}}{t_{a_i}-t_{a_j}}
\prod_{i=1}^r\prod_{j=1}^{\pr} \frac{1}{t_{a_i}-t_{b_j}}
\prod_{a<b}\frac{1}{t_b-q^{-2}t_a}
\\
\quad{} \times \prod_{i=1}^r\delta_{\ep_{a_i},+}
\sum_{\ep_{b_1},\dots ,\ep_{b_{\pr}}=\pm}\quad
\prod_{j=1}^{\pr}\ep_{b_j}
\prod_{i=1}^r\prod_{j=1}^{\pr}(t_{a_i}-q^{-1-\ep_{b_j}}t_{b_j})
\prod_{a<b}(q^{\ep_b}t_b-q^{\ep_a}t_a).
\end{gather*}
Let us calculate the sum in $\{\ep_{b_i}\}$ assuming
$\ep_{a_i}=+$. Using
\begin{gather*}
(t_a-q^{-1-\ep_b}t_b)(t_a-q^{-1+\ep_b}t_b)=(t_a-t_b)(t_a-q^{-2}t_b),
\end{gather*}
and
\begin{gather*}
\prod_{a<b}(q^{\ep_b}t_b-q^{\ep_a}t_a) = \prod_{i,j\atop
a_i>a_j}q(t_{a_i}-t_{a_j}) \prod_{i,j\atop
b_i>b_j}(q^{\ep_{b_i}}t_{b_i}-q^{\ep_{b_j}}t_{b_j})
\\
\phantom{\prod_{a<b}(q^{\ep_b}t_b-q^{\ep_a}t_a)=}{}\times
\prod_{i,j\atop a_i>b_j}(qt_{a_i}-q^{\ep_{b_j}}t_{b_j})
\prod_{i,j\atop b_j>a_i}(q^{\ep_{b_j}}t_{b_j}-qt_{a_i}),
\end{gather*}
we have
\begin{gather}
\sum_{\ep_{b_1},\dots ,\ep_{b_{\pr}}=\pm}\quad
\prod_{j=1}^{\pr}\ep_{b_j} \,\,
\prod_{i=1}^r\prod_{j=1}^{\pr}(t_{a_i}-q^{-1-\ep_{b_j}}t_{b_j})
\prod_{a<b}(q^{\ep_b}t_b-q^{\ep_a}t_a) \non
\\
\quad{}= (-1)^{lr-\frac{1}{2}r(r-1)-\sum_{i=1}^r a_i} q^{r\pr
+\frac{1}{2}r(r-1)}
\prod_{i=1}^r\prod_{j=1}^{\pr}(t_{a_i}-t_{b_j})(t_{a_i}-q^{-2}t_{b_j})
\prod_{i,j\atop a_i>a_j}(t_{a_i}-t_{a_j}) \non
\\
\qquad{} \times \sum_{\ep_{b_1},\dots ,\ep_{b_{\pr}}=\pm}
\prod_{j=1}^{\pr}\ep_{b_j} \prod_{i,j\atop
b_i>b_j}(q^{\ep_{b_i}}t_{b_i}-q^{\ep_{b_j}}t_{b_j}). \label{sum-1}
\end{gather}

\begin{lemma} For $N\geq 1$ we have
\begin{gather}
\sum_{\ep_1,\dots ,\ep_N=\pm}\quad \prod_{j=1}^N\ep_j
\prod_{i>j}(q^{\ep_i}t_i-q^{\ep_j}t_j)=0. \label{eq}
\end{gather}
\end{lemma}

\begin{proof} Let
\begin{gather*}
{\mathbf a}_i(\ep)={}^t(1,q^{\ep}t_i,(q^{\ep}t_i)^2,\dots
,(q^{\ep}t_i)^{N-1}).
\end{gather*}
Then the left hand side of \eqref{eq} is equal to
\begin{gather}
\sum_{\ep_1,\dots ,\ep_N=\pm}\quad \prod_{j=1}^N\ep_j
\det\left({\mathbf a}_1(\ep_1),\dots ,{\mathbf a}_N(\ep_N)\right)
= \det\left(\sum_{\ep_1}\ep_1{\mathbf a}_1(\ep_1),\dots ,
\sum_{\ep_N}\ep_N{\mathbf a}_N(\ep_N)\right). \label{eq-1}
\end{gather}
Since
\begin{gather*}
\sum_{\ep}\ep{\mathbf a}_i(\ep) = {}^t\left( 0,(q-q^{-1})t_i,\dots
,(q^{N-1}-q^{-(N-1)})t_i^{N-1} \right),
\end{gather*}
the right hand side of \eqref{eq-1} is zero.
\end{proof}

By this lemma the right hand side of \eqref{sum-1} becomes zero if
$\pr>0$. Consequently $G^{(\nu)}_{(\mu)}=0$ for $\pr>0$. Suppose
that $\pr=0$. In this case $r=l$, $l_i=i$  $(1\leq i\leq l)$ and
\begin{gather*}
(-q^2)^l(q-q^{-1})^{-l}G_{(-^l)}^{(\nu)}(t,z) \non
\\
\qquad{} =\prod_{a<b}\frac{q(t_b-t_a)}{t_b-q^{-2}t_a}
\sum_{a_1,\dots ,a_l=1 \atop a_i\neq a_j (i\neq j)}^l
\prod_{i=1}^l \left( \frac{t_{a_i}}{z_{k_{i}}-qt_{a_i}}
\prod_{j<k_{i}}\frac{z_j-q^{-1}t_{a_i}}{z_j-qt_{a_i}}
\prod_{i>j}\frac{t_{a_i}-q^{-2}t_{a_j}}{t_{a_i}-t_{a_j}} \right).
\end{gather*}
The theorem easily follows from this.

\section{Concluding remarks}

In this paper we study the solutions of the qKZ equation taking
the value in the tensor product of the two dimensional evaluation
representation of $U_q(\widehat{sl}_2)$. The integral formulae are
derived for the highest to highest matrix elements for certain
intertwining operators by using free f\/ield realizations. The
integrals with respect to $u$ variables corresponding to the
operator $J^{-}(u)$ are calculated and the sum arising from the
expression of $J^{-}(u)$ and the screening operator $S(t)$ is
calculated. The formulae thus obtained coincide with those of
Tarasov and Varchenko. The calculations in this paper can be
extended to the case where the vector space $V^{(1)\otimes n}$ is
replaced by a tensor product of more general representations. It
is an interesting problem to perform similar calculations for
other quantum af\/f\/ine algebras \cite{AOS} and the elliptic
algebras \cite{K2}.

In Tarasov--Varchenko's theory solutions of a qKZ equation are
parametrized by elements of the elliptic hypergeometric space
${\cal F}_{ell}$ while the matrix elements are specif\/ied by
intertwiners. It is an interesting problem to establish a
correspondence between intertwining operators and elements of
${\cal F}_{ell}$. With the results of the present paper one can
begin to study this problem. Study in this direction will provide
a new insight on the space of local f\/ields and correlation
functions of integrable f\/ield theories and solvable lattice
models. The corresponding problem in CFT is studied in~\cite{F}.

\appendix

\pdfbookmark[1]{Appendix. List of OPE's}{appendix}

\section*{Appendix. List of OPE's}
Here we list OPE's which are necessary in this paper. Almost all
of them are taken from the paper \cite{KQS}. Let
\begin{gather*}
C(z)=\frac{(q^{-2}z;p)_\infty}{(q^{2}z;p)_\infty}, \qquad
(z)_\infty=(z;p)_\infty.
\\
S_{\ep_1}(t_1)S_{\ep_2}(t_2)=(q^{-2}t_1)^{4s}q^{\ep_1}
\frac{t_1-q^{\ep_2-\ep_1}t_2}{t_1-q^{-2}t_2}
C(t_2/t_1):S_{\ep_1}(t_1)S_{\ep_2}(t_2):, \qquad
|q^{-2}t_2|<|t_1|,
\\
\phi_{-}(z)S_{\ep}(t)=(q^{k}z)^{-s}
\frac{(qt/z)_\infty}{(q^{-1}t/z)_\infty} :\phi_{-}(z)S_{\ep}(t):,
\qquad |q^{-1}t|<|z|,
\\
J^{-}_{\mu}(u)S_{\ep}(t)=q^{-\mu}
\frac{u-q^{-\mu(k+1)-\ep}t}{u-q^{-\mu(k+2)}t}
:J^{-}_{\mu}(u)S_{\ep}(t):, \qquad |q^{-(k+2)}t|<|u|,
\\
\phi_{-}(z)J^{-}_{\mu}(u)= \frac{z-q^{\mu-2-k}u}{z-q^{-1-k}u}
:\phi_{-}(z)J^{-}_{\mu}(u):, \qquad |u|<|q^{k+3}z| \quad
\mbox{for} \ \ \mu=-,
\\
J^{-}_{\mu}(u)\phi_{-}(z)=q^\mu \frac{u-q^{k+2-\mu}z}{u-q^{k+3}z}
:\phi_{-}(z)J^{-}_{\mu}(u):, \qquad |q^{k+1}z|<|u| \quad
\mbox{for} \ \ \mu=+,
\\
[\phi_{-}(z),J^{-}_{\mu}(u)]_q=
\frac{(1-q^2)u(z-q^{\mu-2-k}u)}{(z-q^{-1-k}u)(u-q^{k+3}z)}
:\phi_{-}(z)J^{-}_{\mu}(u):,
\\
J^{-}_{\mu_1}(u_1)J^{-}_{\mu_2}(u_2)=q^{-\mu_1}
\frac{u_1-q^{\mu_1-\mu_2}u_2}{u_1-q^{-2}u_2}
:J^{-}_{\mu_1}(u_1)J^{-}_{\mu_2}(u_2):, \qquad |q^{-2}u_2|<|u_1|,
\\
\phi_{-}(z_1)\phi_{-}(z_2)=(q^{k}z_1)^s\xi(z_1/z_2)
:\phi_{-}(z_1)\phi_{-}(z_2):, \qquad |pz_2|<|z_1|.
\end{gather*}

\begin{note} After completing the paper we were informed that
H.~Awata, S.~Odake and J.~Shiraishi obtained a similar result to
this paper. The result is reviewed in Shiraishi's PhD
thesis~\cite{S2} in which one statement is a conjecture. However
the complete version containing all proofs for all statements had
not been published after all. We would like to thank J.~Shiraishi
for the kind correspondence.
\end{note}

\subsection*{Acknowledgements}

We would like to thank Hitoshi Konno and Yasuhiko Yamada for
valuable discussions and comments. We are also grateful to Atsushi
Matsuo for useful comments on the manuscript.

\pdfbookmark[1]{References}{ref}
\LastPageEnding

\end{document}